\documentclass[11pt]{amsart}

\usepackage[margin=1.25in]{geometry}

\usepackage{amsmath, amssymb, amsthm, graphicx, enumerate, tikz, float, color}
\usepackage{mathtools, comment}
\usepackage{wasysym}
\usepackage[misc]{ifsym}
\usepackage[colorlinks]{hyperref}

\hypersetup{citecolor=blue}
\usetikzlibrary{matrix,arrows,decorations.pathmorphing}

\newtheorem{theorem}{Theorem}[section]
\newtheorem{lemma}[theorem]{Lemma}

\newtheorem{claim}[theorem]{Claim}
\newtheorem*{thmmain}{Main Theorem}

\theoremstyle{definition}
\newtheorem{definition}[theorem]{Definition}
\newtheorem{example}[theorem]{Example}

\theoremstyle{remark}

\numberwithin{equation}{section}

\makeatletter
\@namedef{subjclassname@2020}{%
  \textup{2020} Mathematics Subject Classification}
\makeatother

\begin{document}

\allowdisplaybreaks

\title[Character degree graphs of solvable groups]{A family of graphs that cannot occur as character degree graphs of solvable groups}

\author[M.W. Bissler]{Mark W. Bissler}
\address{Department of General Education, Western Governors University, Salt Lake City, UT 84107}
\email{mark.bissler@wgu.edu}

\author[J. Laubacher]{Jacob Laubacher}
\address{Department of Mathematics, St. Norbert College, De Pere, WI 54115}
\email{jacob.laubacher@snc.edu}

\author[M.L. Lewis]{Mark L. Lewis}
\address{Department of Mathematical Sciences, Kent State University, Kent, OH 44242}
\email{lewis@math.kent.edu}

\subjclass[2020]{Primary 20C15; Secondary 05C25, 20D10}

\date{\today}

\keywords{Character degree graphs, solvable groups, family of graphs\\\indent\emph{Corresponding author.} Mark W. Bissler \Letter~\href{mailto:mbissle2@outlook.com}{mbissle2@outlook.com}}

\begin{abstract}
We investigate character degree graphs of solvable groups. In particular, we provide general results that can be used to eliminate which degree graphs can occur as solvable groups. Finally, we show a specific family of graphs cannot occur as a character degree for any solvable group.
\end{abstract}

\maketitle

\section{Introduction}

Throughout this paper, $G$ will be a finite solvable group. We will write Irr$(G)$ for the set of irreducible characters of $G$, and cd$(G)=\{\chi(1)\mid \chi\in\text{Irr}(G)\}$. Denote $\rho(G)$ to be the set of primes that divide degrees in cd$(G)$ for the character degrees of $G$. The degree graph of $G$, written $\Delta(G)$, is the graph whose vertex set is $\rho(G)$. Two vertices $p$ and $q$ of $\rho(G)$ are adjacent in $\Delta(G)$ if there exists $a\in\text{cd}(G)$ where $pq$ divides $a$. We identify each vertex of a graph with a prime in $\rho(G)$, and we interchange these throughout this paper to avoid having to specify that for a prime $p\in\rho(G)$ associated with a vertex $v$ of a graph. This type of graph has been studied in a variety of places (see \cite{mjl}, \cite{lewis1}, \cite{lewis2}, \cite{palfy}, \cite{zhang}).

In this paper, we introduce techniques to eliminate a graph from occurring as the degree graph of a solvable group. Using these techniques and induction, we show that a specific family of graphs cannot occur as degree graphs of solvable groups. In \cite{palfy}, P\'{a}lfy showed that if $G$ is a solvable group, then for every three vertices in $\rho(G)$, there is some edge in $\Delta(G)$ incident to two of these vertices. With this in mind, we say that a graph $\Gamma$ satisfies P\'{a}lfy's condition if for every three vertices there is some edge incident to two of them. Note this implies that a graph satisfying P\'{a}lfy's condition and that is disconnected must have two complete connected components. This also forces the graph to have at most diameter three. We rely heavily on the results of \cite{sass} when dealing with a graph or subgraph that arises with diameter three. The results in \cite{sass} are from her dissertation work, and they can also be found in \cite{sassy}.

The following main theorem of the paper will be proven using techniques similar to those found in \cite{lewis1}. In particular, we show that there is an infinite family of graphs that satisfy P\'{a}lfy's condition and do not occur as $\Delta(G)$ for solvable group $G$.
	
\begin{thmmain}
Let $\Gamma$ be a graph satisfying P\'{a}lfy's condition, with $k\geq 5$ vertices. Assume that there exist two vertices $p_1$ and $p_2$ in $\Gamma$, such that $p_1$ and $p_2$ are of degree two, $p_1$ is adjacent to $p_2$, and they share no common neighbor. Then $\Gamma$ is not the prime character degree graph of any solvable group.
\end{thmmain}

We do not require that there is an edge between the two vertices adjacent to $p_1$ and $p_2$, respectively. We will see that this implies that there are two graphs for each value $k$ that satisfy the hypotheses of this theorem.

The results from this paper are used, generalized, and applied by the authors in \cite{bisslaub}, \cite{corey}, \cite{groot}, and \cite{jm}.

\section{General Lemmas}\label{sec2}

We say $\Gamma$ is a subgraph of $\Delta$ if $\Gamma$ consists of a subset of the vertices of $\Delta$ and a subset of the edges of $\Delta$ which are incident to only the vertices of $\Gamma$. If either the set of vertices or set of edges of $\Gamma$ is a proper subset of the set of vertices or set of edges of $\Delta$, then we say $\Gamma$ is a proper subgraph of $\Delta$.

\begin{definition}
A vertex $v$ of a graph $\Gamma$ is \textbf{admissible} if the subgraph of $\Gamma$ obtained by removing the vertex $p$ and all edges incident to $p$ does not occur as $\Delta(G)$ for any solvable group $G$, and none of the subgraphs of $\Gamma$ obtained by removing one or more of the edges incident to $p$ occurs as $\Delta(G)$ for any solvable group $G$. 
\end{definition}
Now, using this definition and the minimality of $|G|$ we are able to show for a prime $p$ associated to an admissible vertex in $\Delta(G)$, that $G=O^p(G)$. 

\begin{lemma}\label{first}
Let G be a solvable group, and suppose $p$ is an admissible vertex of $\Delta(G)$. For every proper normal subgroup $H$ of $G$, suppose that $\Delta(H)$ is a proper subgraph of $\Delta(G)$. Then $O^p(G)=G$.
\end{lemma}
\begin{proof}
Assume not; that is, assume that $O^p(G)<G$. By hypothesis, $\Delta(O^p(G))$ is a proper subgraph of $\Delta(G)$. Since $|G:O^p(G)|$ is a power of $p$, the only vertex of $\Delta(G)$ that could be missing is $p$, and any edges that are missing must be incident to $p$. Therefore, $\Delta(O^p(G))$ is a subgraph of $\Delta(G)$ obtained by deleting $p$ or by deleting one or more of the edges incident to $p$. Since $p$ was assumed to be admissible, this cannot occur. Thus, $O^p(G)=G$.
\end{proof}

Next, we show that for any graph $\Delta$ in which every vertex is admissible, $\Delta$ is not the character degree graph of any solvable group.

\begin{lemma}
If $\Gamma$ is a graph in which every vertex is admissible, then $\Gamma$ is not $\Delta(G)$ for any solvable group $G$.
\end{lemma}
\begin{proof}
Let $G$ be a counterexample with $|G|$ minimal. Since $\Delta(G)$ is not empty, $G>1$, and since $G$ is solvable, there must be a prime $p$ so that $O^p(G)<G$. We have assumed, though, that every vertex of $\Delta(G)$ is admissible, and by Lemma \ref{first} this means by the minimality of $G$, that $O^p(G)=G$ for every prime $p$, which is a contradiction.
\end{proof}

The next lemma provides us with a way to show $G$ does not have a normal Sylow p-subgroup for a specific vertex. We refine a technique from \cite{lewis} for the following result.

\begin{lemma}\label{third}
Let $\Gamma$ be a graph satisfying P\'{a}lfy's condition. Let $q$ be a vertex of $\Gamma$, and denote $\pi$ to be the set of vertices of $\Gamma$ adjacent to $q$, and $\rho$ to be the set of vertices not adjacent to $q$. Assume that $\pi$ is the disjoint union of nonempty sets $\pi_1$ and $\pi_2$, and assume that no vertex in $\pi_1$ is adjacent in $\Gamma$ to any vertex in $\pi_2$. Let $v$ be a vertex in $\pi_2$ adjacent to an admissible vertex $s$ in $\rho$. Furthermore, assume there exists another vertex $w$ in $\rho$ that is not adjacent to $v$. Let $G$ be a solvable group such that $\Delta(G)=\Gamma$, and assume that for every proper normal subgroup $H$ of $G$, $\Delta(H)$ is a proper subgraph of $\Delta(G)$. Then a Sylow $q$-subgroup of $G$ for the prime associated to $q$ is not normal.
\end{lemma}
\begin{proof}
Assume not; that is, assume that $Q$ is a normal nonabelian Sylow $q$-subgroup of $G$. By Lemma 4.5 of \cite{lewis}, we have that $\rho(G/Q)\subseteq\pi$. If $G/Q$ is abelian, then $O^s(G)<G$, a contradiction, as $s$ was assumed to be admissible, and by Lemma \ref{first}, $O^s(G)=G$. If $\rho(G/Q)\subseteq\pi_1$ or $\rho(G/Q)\subseteq\pi_2$, then by Lemma 4.6 of \cite{lewis}, we have that $G'\subseteq QT$, where $T$ is either a Hall $\pi_1$-subgroup or Hall $\pi_2$-subgroup of $G$. In each case, we have again that $O^s(G)<G$, which is a contradiction. Thus, $\rho(G/Q)$ intersects both $\pi_1$ and $\pi_2$ nontrivially. By Lemma 4.7 of \cite{lewis}, we have two possibilities. The first possibility is that $G/Q$ has a central Hall $\rho$-subgroup, which implies that $O^s(G)<G$, which is a contradiction. So we have the second conclusion of Lemma 4.7, that is, there is a prime $r\in\rho$ such that $G/Q$ has a noncentral Sylow $r$-subgroup, and every other prime $t\in\rho$ is adjacent in $\Delta(G)$ to every prime in $\rho(G/Q)$. If $r\neq w$, we have that $w$ is adjacent in $\Delta(G)$ to every prime in $\rho(G/Q)$, which means that $w$ is adjacent to $v$, a contradiction, as we assumed they were not adjacent. Now when $r=w$, we know that a Sylow $s$-subgroup of $G/Q$ has a normal Hall $s$-complement. This implies that $O^s(G)<G$, and this is a contradiction, as $s$ was assumed to be admissible. Thus, $G$ does not have a normal Sylow $q$-subgroup.
\end{proof}
	
We now make a refinement on the definition of an admissible vertex.

\begin{definition}
We say a vertex $p$ of a graph $\Gamma$ is \textbf{strongly admissible} if $p$ is admissible, and none of the subgraphs of $\Gamma$ obtained by removing $p$, the edges incident to $p$, and one or more of the edges between two adjacent vertices of $p$ occurs as $\Delta (G)$ for some solvable group $G$.
\end{definition}

\begin{lemma}\label{second} Let G be a solvable group, and assume that $p$ is a prime whose vertex is a strongly admissible vertex of $\Delta(G)$. For every proper normal subgroup $H$ of $G$, suppose that $\Delta(G/H)$ is a proper subgraph of $\Delta(G)$. Then a Sylow $p$-subgroup of $G$ is not normal.
\end{lemma}
\begin{proof}
Assume not; that is, suppose $P$ is a normal nonabelian Sylow $p$-subgroup of $G$. By Lemma 3 of \cite{lewis1} we know that $\rho(G/P')=\rho(G)\setminus{\{p\}}$. Let $\{q,r\}$ be an edge of $\Delta(G)$ with $q$ and $r$ distinct from $p$. If $\{q,r\}$ is not an edge of $\Delta(G/P')$, then we will show that $q$ and $r$ are adjacent to $p$. Since $\{q,r\}$ is an edge of $\Delta(G)$, there exists $\chi\in\mbox{Irr}(G)$ with $qr\mid\chi(1)$. We know that $P'\nless\ker\chi$, which implies that $p\mid\chi(1)$. Thus, $q$ and $r$ are adjacent to $p$. The graph $\Delta (G/P')$ can be obtained from the graph $\Delta (G)$ by removing $p$ and the edges incident to $p$, and perhaps also by removing one or more edges that are incident to vertices that are adjacent to $p$. Since $p$ was assumed to be strongly admissible, we know the resulting subgraphs cannot occur as the prime character degree graph of a solvable group. Therefore, no such graph can be the prime character degree graph of a quotient of $G$; thus, a Sylow $p$-subgroup of $G$ is not normal in $G$.
\end{proof}

By \cite{lewis2}, we know the exact structure of a solvable group $G$, whose character degree graph has two connected components. Furthermore in \cite{lewis2}, it is shown that if $G$ is a solvable group, then $\Delta(G)$ has two connected components if and only if $G$ is one of the groups from Examples 2.1-2.6. The graphs we observe in this paper satisfy either Example 2.4 or Example 2.6 of \cite{lewis2}. That is:

\begin{example}[Example 2.4 from \cite{lewis2}]
The group $G$ is the semi-direct product of a subgroup $H$ acting on an elementary abelian $p$-group $V$ for some prime $p$. Let $Z=C_H(V)$ and $K$ be the Fitting subgroup of $H$. Write $m=|H:K|>1$, and $|V|=q^m$ where $q$ is a $p$-power. We have $Z\subseteq Z(H)$, $K/Z$ is abelian, $K$ acts irreducibly on $V$, $m$ and $|K:Z|$ are relatively prime, and that $(q^{m}-1)/(q-1)$ divides $|K:Z|$.
\end{example}

\begin{example}[Example 2.6 from \cite{lewis2}]
The group $G$ is the semi-direct product of an abelian group $D$ acting coprimely on a group $T$ so that $[T,D]$ is a Frobenius group. The Frobenius kernel is $A=T'=[T,D]'$, $A$ is a non-abelian $p$-group for some prime $p$, and a Frobenius complement is $B$ with $[B,D]\subseteq B$. Every character in $\mbox{ Irr}(T\mid A')$ is invariant under the action of $D$ and $A/A'$ is irreducible under the action of $B$. If $m=|D:C_D(A)|$, then $|A:A'|=q^m$ where $q$ is a $p$-power, and $(q^m-1)/(q-1)$ divides $|B|$.
\end{example}

The next lemma deals with a specific size of connected components within a degree graph. We manipulate Lemma 5.1 of \cite{lewis} for the following result using the Zsigmondy prime theorem. That is, let $a>1$ and $n$ be positive integers. Then there exists a Zsigmondy prime divisor for $a^n-1$ unless:
\begin{enumerate}[(1)]
\item $n=2$ and $a=2^k-1$ for some $k\in\mathbb{N}$, or
\item $n=6$ and $a=2$.
\end{enumerate}
For reference on the Zsigmondy prime theorem, we suggest Theorem 6.2 of \cite{wolf}.

\begin{lemma}\label{fourth}
Let G be a solvable group and $\Delta(G)$ be its prime character degree graph. Let $N\triangleleft G$, and assume that $\rho(G/N)=\rho(G)$. Furthermore, assume that $G/N$ satisfies the hypotheses of Example 2.4 and has the two connected components $\pi=\{p_1,p_2\}$ and $\rho=\rho(G)\setminus\pi$. If $p_1$ and $p_2$ have degree two in $\Delta(G)$, and $p_1$ and $p_2$ share no common neighbor aside from each other, then $N=1$. 
\end{lemma}
\begin{proof}
Assume not; that is, assume that $N>1$. Let $F/N$ and $E/F$ be the Fitting subgroups of $G/N$ and $G/F$ respectively. Since we are supposing that $G/N$ satisfies the hypotheses of Example 2.4 from \cite{lewis2}, we know that $|G:E|=p_1^{\alpha_1}p_2^{\alpha_2}$ and that $(|G:E|,|E:F|)=1$. Also, following Lemma 5.2 from \cite{lewis}, there exist primes $q\in\rho$ and $s\in\pi$ such that $(q^{|G:E|}-1)/(q^{\frac{|G:E|}{s}}-1)$ divides $|E:F|$. We now want to show that $(q^{|G:E|}-1)/(q^{\frac{|G:E|}{s}}-1)$ is divisible by at least two primes. We know that $$q^{|G:E|}-1=\prod_{d\mid |G:E|}\Phi_d(q), \mbox{ and}$$  $$q^{\frac{|G:E|}{s}}-1=\prod_{d\mid \frac{|G:E|}{s}}\Phi_d(q), \mbox{ thus}$$  $$(q^{|G:E|}-1)/(q^{\frac{|G:E|}{s}}-1)=\prod_{\substack{d\mid |G:E|} \\, {d\nmid \frac{|G:E|}{s}}} \Phi_d(q).$$ Now we wish to apply the Zsigmondy prime theorem to this quotient. In order to do this, we must check that the exceptions do not occur.

Suppose that $2$ divides $|G:E|$ and $q$ is odd. It follows that $2$ divides $q+1=(q^2-1)/(q-1)$. Now, we also know that $(q^2-1)/(q-1)$ divides $(q^{|G:E|}-1)/(q^{\frac{|G:E|}{s}}-1)$, which in turn divides $|E:F|$. However, this means that $2$ divides both $|G:E|$ and $|E:F|$, a contradiction, as we know these two indices are relatively prime. Thus, either $|G:E|$ is odd, or $q$ is even. 

If $q=2$ and $6$ divides $|G:E|$, then $3=(2^2-1)/(2-1)$ is a divisor of $(2^{|G:E|}-1)/(2^{\frac{|G:E|}{s}}-1)$, which divides $|E:F|$. This is a contradiction again, as $(|G:E|,|E:F|)=1$. Now, by the Zsigmondy prime theorem, we have that $(q^{|G:E|}-1)/(q^{\frac{|G:E|}{s}}-1)$ has at least two prime divisors. By our assumptions, since $p_1$ is adjacent to only $p_2$ and one vertex in $\rho$ which is not adjacent to $p_2$, we can apply Lemma 5.2 of \cite{lewis} to see that $p_1$ must be adjacent to all the primes that divide $(q^{|G:E|}-1)/(q^{\frac{|G:E|}{s}}-1)$. If $p_2$ is one of these divisors, then $p_1$ and $p_2$ would have a common neighbor, a contradiction as we assumed they shared no common neighbors. If neither of the two common divisors is $p_2$, then this implies $p_1$ would have degree at least $3$,  which contradicts the fact that $p_1$ was assumed to have degree $2$ in $\Delta(G)$.
\end{proof}

The last lemma we provide in this section will be the final tool we use to show a graph cannot occur for any solvable group $G$. We address some notation before we begin. For  $\theta\in$\text{Irr}$(N)$ with $N\triangleleft G$, we note the standard notation to define Irr$(G|\theta)$ is to be the set of characters in Irr$(G)$ that are constituents of $\theta^G$. Following this notation, we define cd$(G|\theta)=\{\chi(1)\mid\chi\in\text{Irr}(G\vert\theta)\}$. In the same manner, we define Irr$(G|N)$ to be the union of the sets Irr$(G\vert\theta)$, where $\theta$ runs through all the nonprincipal characters in Irr$(N)$. Our interest is with the set cd$(G|N)=\{\chi(1)\vert\chi\in\text{Irr}(G\vert N)\}$.

\begin{lemma}\label{final}
Let $\Gamma$ be a graph satisfying P\'{a}lfy's condition with $n\geq 5$ vertices and no complete vertices. Also, assume there exist vertices $p_1$ and $p_2$ of $\Gamma$ such that $p_1$ is adjacent to an admissible vertex $q_1$ and $p_2$ is not adjacent to $q_1$, and $p_1$ is not adjacent to $q_2$, another admissible vertex. Let $G$ be a solvable group and suppose for all proper normal subgroups $N$ of $G$ we have that $\Delta(N)$ and $\Delta(G/N)$ are proper subgraphs of $\Gamma$. Let $F$ be the Fitting subgroup of $G$ and suppose $F$ is a minimal normal subgroup of $G$. Then $\Gamma$ is not the prime character degree graph of $G$.
\end{lemma}
\begin{proof}
Since $F$ was assumed to be a minimal normal subgroup of $G$, we have that $\Phi(G)=1$, and by Lemma III 4.4 of \cite{huppert}, there is a subgroup $H$ so that $G=HF$, and $H\cap F=1$. Denote $E$ to be the Fitting subgroup of $H$. Let $p$ be a prime divisor of $|E|$. We know via Lemma 2.10 of \cite{lewis3} that every degree in $\mbox{cd}(G|F)$ is divisible by all the prime divisors of $|E|$. By assumption, since there are no complete vertices, we know there is a prime $q\in\rho(G)$ that is not adjacent to $p$ in $\Delta(G)$. Consider a non-principal character $\lambda\in\mbox{Irr}(F)$. We know that every degree in $\mbox{cd}(G|\lambda)$ is divisible by $p$, so $C_H(\lambda)$ contains a Sylow $q$-subgroup of $H$ as a normal subgroup. Using Lemma 1 of \cite{lewis4}, either $(1)$  $ H\cong \mbox{SL}_2(3)$ or $\mbox{Gl}_2(3)$, or $(2)$ $H$ has a normal abelian subgroup that acts irreducibly on $F$. We know that $\rho(G)=\pi(|H|)$. If $(1)$ occurs, then $\pi(|H|)=\{2,3\}$, which is not allowed since, $|\rho(G)|\geq 5$. Thus, $(2)$ must occur.

In case $(2)$, $H/E$ is abelian, and $|H:E|\in \mbox{cd}(G)$. There is a degree in $\mbox{cd}(G)$ divisible by all the prime divisors of $|E|$. Any prime in $\pi(|H:E|)\cap\pi(|E|)$ would be adjacent in $\Delta(G)$ to all of the primes in $\rho(G)$. Since $\Delta(G)$ has no such vertex, $|H:E|$ and $|E|$ are relatively prime.

If $q_1$ divides $|H:E|$, then $O^{q_1}(H)<H$, and thus $O^{q_1}(G)<G$. Since $q_1$ was assumed admissible, this cannot happen, so $q_1$ does not divide $|H:E|$, and $q_1$ divides $|E|$. We have that $p_2$ divides $|H:E|$, and that $q_2$ divides $|E|$. This yields $\pi(|H:E|)=\{p_1,p_2\}$, and $\pi(|E|)=\rho(G)\setminus \{p_1,p_2\}$.

Following a generalized argument of Claim 3 from \cite{lewis}, we next let $Q_1$ be a Sylow $q_1$-subgroup and $Q_2$ be a Sylow $q_2$-subgroup of $E$. Let $A$ be a Hall $q_1$-complement for $H$, so $H=Q_1A$ and $Q_1\cap A=1$. We can find a character $\chi\in \mbox{Irr}(G)$ with $p_1q_1$ dividing $\chi(1)$. We know that $q_2$ does not divide $\chi(1)$, so $\chi\notin \mbox{Irr}(G|F)$. It follows that $\chi\in\mbox{Irr}(G/F)$ and $\chi_H$ is irreducible. Let $\theta$ be an irreducible constituent of $\chi_{Q_1}$. Now $q_1$ divides $\theta(1)$, so $Q_1$ is not abelian. The stabilizer of $\theta$ in $H$ is $Q_1C_A(\theta)$. Observe that $Q_2\subseteq C_A(\theta)$. Using the usual arguments, $p_2\notin\rho(C_A(\theta))$, and $C_A(\theta)$ contains a Sylow $p_2$-subgroup $P_2$ of $H$ as a normal abelian subgroup. We see that $P_2$ centralizes $Q_2$, and $P_2E=Q_1TP_2\times Q_2$, where $T$ is the Hall $\{q_1, q_2\}$-subgroup of $E$. Since $H/E$ is abelian, $P_2E$ is normal in $H$. Let $K=Q_1TP_2$, and note that $FK$ is normal in $G$. Furthermore, $\Delta(FK)$ has two connected components, and thus $FK$ is one of the examples of two connected components mentioned prior (Examples 2.1-2.6 from \cite{lewis2}), but this is a contradiction, as the only examples where $EF/F$ is not abelian are Examples 2.2 and 2.3, and in both of those cases $\rho(G)=\{2,3\}$. This is the final contradiction, and the theorem is proved.
\end{proof}

\section{Infinite Family}

We now use techniques from Section \ref{sec2} to show that an infinite family of graphs cannot occur for any solvable group $G$, which is detailed in the Main Theorem again stated below.

\begin{thmmain}
Let $\Gamma$ be a graph satisfying P\'{a}lfy's condition, with $k\geq 5$ vertices. Assume that there exist two vertices $p_1$ and $p_2$ in $\Gamma$, such that $p_1$ and $p_2$ are of degree two, $p_1$ is adjacent to $p_2$, and they share no common neighbor. Then $\Gamma$ is not the prime character degree graph of any solvable group.
\end{thmmain}

We do not require that there is an edge between the two vertices adjacent to $p_1$ and $p_2$, respectively. We will see that this implies that there are two graphs for each value of $k$ that satisfy the hypotheses of this theorem. We provide an example of this in Figure \ref{fig1} to illustrate the two types of graphs arising.

\begin{figure}[h]
\centering
\begin{tikzpicture}[scale=2]
\node (1a) at (0,1) {$p_1$};
\node (1b) at (0,0) {$p_2$};
\node (1c) at (1.5,1) {$q_1$};
\node (1d) at (1.5,0) {$q_2$};
\node (1e) at (2,.5) {$\bullet$};
\node (1f) at (2.5,.5) {$\bullet$};
\path[font=\small,>=angle 90]
(1a) edge node [right] {$ $} (1b)
(1c) edge node [right] {$ $} (1d)
(1c) edge node [right] {$ $} (1e)
(1c) edge node [right] {$ $} (1f)
(1d) edge node [right] {$ $} (1e)
(1d) edge node [right] {$ $} (1f)
(1e) edge node [right] {$ $} (1f)
(1a) edge node [right] {$ $} (1c)
(1b) edge node [right] {$ $} (1d);
\node (2a) at (4,1) {$p_1$};
\node (2b) at (4,0) {$p_2$};
\node (2c) at (5.5,1) {$q_1$};
\node (2d) at (5.5,0) {$q_2$};
\node (2e) at (6,.5) {$\bullet$};
\node (2f) at (6.5,.5) {$\bullet$};
\path[font=\small,>=angle 90]
(2a) edge node [right] {$ $} (2b)
(2c) edge node [right] {$ $} (2e)
(2c) edge node [right] {$ $} (2f)
(2d) edge node [right] {$ $} (2e)
(2d) edge node [right] {$ $} (2f)
(2e) edge node [right] {$ $} (2f)
(2a) edge node [right] {$ $} (2c)
(2b) edge node [right] {$ $} (2d);
\node (3a) at (0,-.5) {$p_1$};
\node (3b) at (0,-1.5) {$p_2$};
\node (3c) at (1.5,-1) {$\bullet$};
\node (3d) at (2,-.5) {$q_1$};
\node (3e) at (2,-1.5) {$q_2$};
\node (3f) at (2.5,-.7) {$\bullet$};
\node (3g) at (2.5,-1.3) {$\bullet$};
\path[font=\small,>=angle 90]
(3a) edge node [right] {$ $} (3b)
(3c) edge node [right] {$ $} (3d)
(3c) edge node [right] {$ $} (3e)
(3c) edge node [right] {$ $} (3f)
(3c) edge node [right] {$ $} (3g)
(3d) edge node [right] {$ $} (3e)
(3d) edge node [right] {$ $} (3f)
(3d) edge node [right] {$ $} (3g)
(3e) edge node [right] {$ $} (3f)
(3e) edge node [right] {$ $} (3g)
(3f) edge node [right] {$ $} (3g)
(3a) edge node [right] {$ $} (3d)
(3b) edge node [right] {$ $} (3e);
\node (4a) at (4,-.5) {$p_1$};
\node (4b) at (4,-1.5) {$p_2$};
\node (4c) at (5.5,-1) {$\bullet$};
\node (4d) at (6,-.5) {$q_1$};
\node (4e) at (6,-1.5) {$q_2$};
\node (4f) at (6.5,-.7) {$\bullet$};
\node (4g) at (6.5,-1.3) {$\bullet$};
\path[font=\small,>=angle 90]
(4a) edge node [right] {$ $} (4b)
(4c) edge node [right] {$ $} (4d)
(4c) edge node [right] {$ $} (4e)
(4c) edge node [right] {$ $} (4f)
(4c) edge node [right] {$ $} (4g)
(4d) edge node [right] {$ $} (4f)
(4d) edge node [right] {$ $} (4g)
(4e) edge node [right] {$ $} (4f)
(4e) edge node [right] {$ $} (4g)
(4f) edge node [right] {$ $} (4g)
(4a) edge node [right] {$ $} (4d)
(4b) edge node [right] {$ $} (4e);
\end{tikzpicture}
\caption{Example of graphs from the Main Theorem}
\label{fig1}
\end{figure}

\begin{proof}
Let $G$ be a counterexample with $|G|$ minimal. We proceed by induction on $|\rho(G)|$ with our original assumptions. Since two graphs arise under our hypotheses, whether $q_1$ and $q_2$ are adjacent in $\Gamma$, we will show simultaneously that both cannot occur. Note that when we consider the case where $q_1$ and $q_2$ are adjacent in $\Gamma$, we will have to assume that the previous graph without them adjacent cannot occur. Note that the graphs with five vertices satisfying our hypotheses were shown not to occur in Corollary 3.3 and Theorem 6.1 of \cite{lewis}.

First, for primes in $\rho(G)$, the vertex adjacent to $p_1$ that is not $p_2$ will be labelled $q_1$, and similarly label $q_2$ the vertex adjacent to $p_2$ that is not $p_1$. Let $r$ and $s$ be distinct vertices in $\rho(G)\setminus\{p_1,p_2,q_1,q_2\}$. By P\'{a}lfy's condition, $r$ and $s$ must be adjacent to every vertex in $\Gamma$ except $p_1$ and $p_2$. We first show which vertices of our graph are admissible.

\begin{claim}
Every vertex in $\rho(G)$ is strongly admissible except $p_1$ and $p_2$.
\end{claim}
\begin{proof}
First, we show that the vertices in $\rho(G)\setminus\{p_1,p_2,q_1,q_2\}$ are strongly admissible, which we denote arbitrarily by $r$. If $r$ loses the edge with $q_1$ or $q_2$,  we violate P\'{a}lfy's condition with $p_2, q_1, r$ and $p_1, q_2, r$, respectively. If the edge between $r$ and $s$ is removed, again P\'{a}lfy's condition is violated with $p_1$, $r,$ and $s$. Now, by our inductive assumption, if $r$ is removed from $\Gamma$, we arrive at a subgraph with one fewer vertex, which cannot occur. Thus, $r$ is admissible and, therefore, any vertex in $\rho(G)\setminus\{p_1,p_2,q_1,q_2\}$ is also admissible for the same reason. If $q_1$ loses the edge with $p_1$, we arrive at a graph with diameter three, which is not possible by Corollary 5.5 of \cite{sass}; similarly for $q_2$ and $p_2$. Next, consider the case when $q_1$ and $q_2$ are adjacent. If the edge is lost between these two vertices, we will arrive at the first graph, which we assume does not occur when these two vertices are adjacent. This shows that any vertex in $\rho(G)\setminus\{p_1,p_2,q_1,q_2\}$ is strongly admissible.

Now if $q_1$ is removed from $\Gamma$, we arrive at a subgraph with diameter three, which again is not possible by Corollary 5.5 of \cite{sass}. The only other edges that could be possibly lost in this subgraph would be between two vertices that are adjacent to $q_1$. However, in every one of these cases, we arrive at a connected or even disconnected graph that violates P\'{a}lfy's condition. Thus, $q_1$ is strongly admissible, and by symmetry of the graph, $q_2$ is strongly admissible.
\end{proof}

We have shown that all vertices, aside from $p_1$ and $p_2$, are strongly admissible. We know by Lemma \ref{second} that the corresponding Sylow subgroups of $G$ are not normal. Next, we show that $G$ does not have a normal Sylow $p_1$- or $p_2$-subgroup.

\begin{claim}
The group $G$ does not have a normal Sylow $p_1$-subgroup or Sylow $p_2$-subgroup.
\end{claim}
\begin{proof}
We will show that $p_1$ satisfies the hypotheses of Lemma \ref{third}. Note we have $\rho=\rho(G)\setminus\{p_1,p_2,q_1\}$ and $\pi=\{p_2,q_1\}$. Since $p_2$ and $q_1$ are not adjacent, we have that $\pi=\pi_1\cup\pi_2$, where $\pi_1=\{q_1\}$ and  $\pi_2=\{p_2\}$. We have previously shown that $q_1$ and $q_2$ are strongly admissible, and note that $p_2$ is adjacent to $q_2$. The last hypothesis we verify is that there exists another vertex in $\rho$ that is not adjacent to $p_2$. Since $\Gamma$ was assumed to have at least five vertices, we know there exists another vertex in $\Gamma$ that is not adjacent to $p_2$. Thus, we satisfy the hypotheses of Lemma \ref{third}, and therefore $G$ does not have a normal Sylow $p_1$-subgroup. By symmetry of $\Gamma$, a similar argument works for $p_2$.
\end{proof}

We have shown that $G$ has no normal nonabelian Sylow subgroups. Let $F$ be the Fitting subgroup of $G$. We note that $\rho(G)=\pi(|G:F|)$, and thus $\rho(G)=\rho(G/\Phi(G))$, where $\Phi(G)$ is the Frattini subgroup of $G$. We next work to show $\Phi(G)=1$. The following claim is under the assumption that $q_1$ and $q_2$ are not adjacent in $\Gamma$. Afterwards, we will provide another claim for the case when they are adjacent.

\begin{claim}
Suppose $H$ is a group where $\rho(H)=\rho(G)$, $\Delta(H)$ is a subgraph of $\Delta(G)$, and $|H|<|G\mid.$ Then $\Delta(H)=\Delta(G)$.
\end{claim}
\begin{proof}
Assume not, that is, $\Delta(H)\neq\Delta(G)$. Applying P\'{a}lfy's condition to $p_1$, $p_2$, and $r$ in $\Delta(H)$, we see that $p_1$ and $p_2$ must be adjacent in $\Delta(H)$. If $\Delta(H)$ is disconnected, P\'{a}lfy's condition forces each component to be a complete graph, and since our graph does not have two complete components, we may assume that $\Delta(H)$ is connected. If $p_1$ is not adjacent to $q_1$ in $\Delta(H)$, we have a graph that is diameter three; by Corollary 5.5 of \cite{sass}, this is not possible. Similarly, if $p_2$ is not adjacent to $q_2$, we arrive at a graph with diameter three, which is not possible for the same reason as before. If $q_1$ is not adjacent to $r$, then P\'{a}lfy's condition is violated for $p_2$, $q_1$, and $r$. If $q_2$ is not adjacent to $r$, then we again violate P\'{a}lfy's condition with $p_1$, $q_2$, and $r$. Lastly, if we lose an edge between $r$ and $s$, we also violate P\'{a}lfy's conditon with $r$, $s$, and $p_1$. Since $\Delta(H)$ cannot lose any edge, we must have that $\Delta(H)=\Delta(G)$.
\end{proof}

Now we address the case when $q_1$ and $q_2$ are adjacent in $\Gamma$.

\begin{claim}\label{cl614}
Suppose $H$ is a group where $\rho(H)=\rho(G)$,  $\Delta(H)$ is a subgraph of $\Delta(G)$, and $|H|<|G|.$ Then $\Delta(H)$ either has connected components,
 $\{p_1,p_2\}$ and $\rho(G)\setminus\{p_1,p_2\}$, or $\Delta(H)=\Delta(G)$.
\end{claim}
\begin{proof}
Applying P\'{a}lfy's condition to $p_1$, $p_2$, and $r$ in $\Delta(H)$, we see that $p_1$ and $p_2$ must be adjacent in $\Delta(H)$. If $\Delta(H)$ is disconnected, P\'{a}lfy's condition forces each component to be a complete graph, and so the two connected components must be $\{p_1,p_2\}$ and $\rho(G)\setminus\{p_1,p_2\}$. So now we assume that $\Delta(H)$ is connected, and we show that it is in fact $\Delta(G)$. If $q_1$ and $q_2$ are not adjacent, we arrive at the graph already shown not to occur. If $p_1$ is not adjacent to $q_1$ in $\Delta(H)$, we have a graph that is diameter three, and by Corollary 5.5 of \cite{sass}, this is not possible. Similarly, if $p_2$ is not adjacent to $q_2$, we arrive at a graph with diameter three, which is not possible for the same reason. If $q_1$ is not adjacent to $r$, then P\'{a}lfy's condition is violated for $p_2$, $q_1$, and $r$. If $q_2$ is not adjacent to $r$, then we again violate P\'{a}lfy's condition with $p_1$, $q_2$, and $r$. Finally, if we lose an edge between $r$ and $s$, we also violate P\'{a}lfy's condition with $r$, $s$, and $p_1$. Since $\Delta(H)$ cannot lose any edge, we must have that $\Delta(H)=\Delta(G)$.
\end{proof}

We now show that if $M\triangleleft G$ with $\rho(G/M)=\rho(G)$, then $M=1$. If $q_1$ is not adjacent to $q_2$, we have this immediately from Claim \ref{cl614}. If $q_1$ is adjacent to $q_2$, we have two cases to consider, $\Delta(G/M)$ is connected and is $\Delta(G)$, or $\Delta(G/M)$ has exactly two connected components. If $\Delta(G/M)$ is connected and is $\Delta(G)$, we are done, as this forces $M=1$. In the disconnected case, we see that this subgraph satisfies Example 2.4 or Example 2.6. First, assume $\Delta(G/M)$ satisfies Example 2.4. Then by Lemma \ref{fourth} we have the desired result, that is, $M=1$. If $\Delta(G/M)$ satisfies Example 2.6, we know that $G/M$ has a normal Sylow $p$-subgroup $P/M$. We also know that $G/P'$ satisfies Example 2.4 by Lemma 3.6 of \cite{lewis2}. Thus, we can again apply Lemma \ref{fourth} to get the desired result that $M=1$.

Since $\rho(G/\Phi(G))=\rho(G)$, we now have that $\Phi(G)=1$. Applying Lemma III 4.4 of \cite{huppert}, there is a subgroup $H$ so that $G=HF$ and $H\cap F=1$. Let $E$ be the Fitting subgroup of $H$. Next we show that the Fitting subgroup of $G$ is minimal normal.

\begin{claim}\label{cl615}
The Fitting subgroup $F$ of $G$ is a minimal normal subgroup.
\end{claim}
\begin{proof}
Suppose that there is a normal subgroup $N$ of $G$ so that $1<N<F$. By Theorem III 4.5 of \cite{huppert}, there is a normal subgroup $M$ of $G$ so that $F=N\times M$. Since $N>1$ and $M>1$, we have $\rho(G/N)\subset\rho(G)$ and $\rho(G/M)\subset\rho(G)$. For any prime $p\in\rho(G)\setminus\rho(G/N)$, we know that $G/N$ has a normal abelian Sylow $p$-subgroup. The class of finite groups with an abelian and normal Sylow $p$-subgroup is a formation, so $p$ must lie in $\rho(G/M)$. Thus, $\rho(G)=\rho(G/N)\cup\rho(G/M)$.

If $p\in\rho(G)\setminus\rho(G/N)$, then $p$ is not in $\rho(G/F)=\rho(H)$; so $E$ must then contain the Sylow $p$-subgroup of $H$. Since $p\in\rho(G)$, it follows that $p$ divides $|H|$, and so $p$ will divide $|E|$. Recall that $\mbox{cd}(G)$ contains a degree divisible by all the prime divisors of $|EF:F|=|E|$. We conclude that $\rho(G)\setminus(\rho(G/M)\cap \rho(G/N))$ lies in a complete subgraph of $\Delta(G)$. Therefore, $\rho(G)\setminus(\rho(G/M)\cap \rho(G/N))$ lies in the subsets: (1) $\{p_1,p_2\}$, (2) $\{p_1,q_1\}$, (3) $\{p_2,q_2\},$ or (4) $\rho(G)\setminus\{p_1,p_2\}$.

Suppose that $(1)$ occurs. This implies that $E$ contains a Hall $\{p_1,p_2\}$-subgroup of $H$. Since $\mbox{cd}(G)$ has a degree divisible by all the primes dividing $|E|$, we see that $|E|$ is divisible by no other primes, and $E$ is the Hall $\{p_1,p_2\}$-subgroup of $H$. We can find a character $\chi\in\mbox{Irr}(G)$ with $p_1p_2$ dividing $\chi(1)$. Let $\theta$ be an irreducible constituent of $\chi_{_{FE}}$. Now, $\chi(1)/\theta(1)$ divides $|G:FE|$ and $\chi(1)$ is relatively prime to $|G:FE|$. We determine that $\chi_{_{FE}}=\theta$. Since $p_1$ and $p_2$ divide $\theta(1)$, and the only possible prime divisors of $a\in\mbox{cd}(G/FE)$ are $\rho(G)\setminus\{p_1,p_2\}$, we conclude via Gallagher's theorem that $\mbox{cd}(G/FE)=\{1\}$ and $G/FE$ is abelian. This implies that $O^{q_1}(G)<G$, which is a contradiction, as we have shown $q_1$ is strongly admissible, and by Lemma \ref{first}, $O^{r}(G)=G$. Thus, $(1)$ cannot occur.

Suppose that $(2)$ occurs. Then we have that $\rho(G)\setminus\{p_1,q_1\}\subseteq\rho(G/N)\cap\rho(G/M)$. We now consider the possible cases for $\rho(G/N)$ and $\rho(G/M)$. Since  $\rho(G)=\rho(G/N)\cup\rho(G/M)$, we must have that $p_1$ is in $\rho(G/N)$ or $\rho(G/M)$. Assume that $p_1\in\rho(G/N)$.  Note that the connected graph with this vertex set would have diameter three from the vertex $p_1$ to $r$, and by Corollary 5.5 of \cite{sass}, we know that this is not possible. Thus, the only subgraph to arise with these vertices will have connected components $\pi=\{p_1, p_2\}$ and $\rho=\rho(G)\setminus\{p_1, p_2, q_1\}$. By Theorem 5.5 of \cite{lewis2}, $G/N$ has a central Sylow $q_1$-subgroup. However, this would imply $O^{q_1}(G)<G$, a contradiction as we have shown $q_1$ is strongly admissible. Thus, (2) cannot occur. By symmetry of the graph, we have also that (3) cannot happen.

Finally, suppose $(4)$ occurs. Then, $\{p_1,p_2\}\subseteq \rho(G/N)\cap\rho(G/M)$. Note that $|\rho(G/M)|\geq 4$ or $|\rho(G/N)|\geq 4$, as we assumed the original graph has at least five vertices. Assume that $|\rho(G/N)|=4$. Let $r\in\rho(G)\setminus\{p_1,p_2,q_1,q_2\}$. If $q_1$ and $r$ are contained in $\rho(G/N)$, then we know by Theorem 1 of \cite{palfy} and Theorem 5 of \cite{zhang} that the graphs arising from this vertex set  cannot occur. Now assume that $q_1$ and $q_2$ are contained in $\rho(G/N)$. Since $\rho(G)=\rho(G/N)\cup\rho(G/M)$, we know $\{p_1,p_2,r\}\subseteq\rho(G/M)$. Then $\Delta(G/M)$ has two connected components, $\{p_1,p_2\}$ and $\{r\}$. Now by Theorem 5.5 of \cite{lewis2}, $G/M$ has either a central Sylow $q_1$-subgroup or a central Sylow $q_2$-subgroup. This would imply either $O^{q_1}(G)<G$ or $O^{q_2}(G)<G$, a contradiction, as we have shown $q_1$ and $q_2$ are both strongly admissible. Thus, $|\rho(G/N)|>4$. If $q_1$ and $q_2$ are both not contained in $\rho(G/N)$, we arrive at the same contradiction via Theorem 5.5 of \cite{lewis2}, with a Sylow $q_1$-subgroup or Sylow $q_2$-subgroup being central in $G/N$. Now assume that $q_1$ is contained in $\rho(G/N)$ and $q_2$ is not contained in $\rho(G/N)$. The connected graph cannot occur by Corollary 5.5 of \cite{sass}, so we have the disconnected graph with two connected components. Again, we apply Theorem 5.5 of \cite{lewis2} to see that a Sylow $q_2$-subgroup is central in $G/N$. This implies that $O^{q_2}(G)<G$, a contradiction, since $q_2$ was shown to be strongly admissible, which implies $O^{p_2}(G)=G$. Note that a similar argument works for the case when $q_2$ is contained in $\rho(G/N)$ and $q_1$ is not contained in $\rho(G/N)$. Now we may assume that both $q_1$ and $q_2$ are contained in $\rho(G/N)$. The connected graph will satisfy our inductive hypothesis and thus cannot occur. The disconnected graph will have the two connected components $\{p_1,p_2\}$ and $\rho(G)\setminus\{p_1,p_2,r\}$. Applying Theorem 5.5 of \cite{lewis2} one last time we have that a Sylow $r$-subgroup is central in $G/N$, which is a contradiction again, since every vertex in $\rho(G)\setminus\{p_1,p_2\}$ was shown to be strongly admissible. This is the final case, so the claim is proved. 
\end{proof}

Now we will show that the hypotheses of Lemma \ref{final} are satisfied to complete the proof of the Main Theorem. First we note that by construction that there are no complete vertices (see Figure \ref{fig1}). Next, we have that $p_1$ is adjacent to $p_2$, and $p_1$ is adjacent to $q_1$, an admissible vertex not adjacent to $p_2$. Also, there exists an admissible vertex $q_2$ that is not adjacent to $p_1$, and from the previous claim, we know that $F$ is minimal normal. Thus, the hypotheses of Lemma \ref{final} are now satisfied, and $\Gamma$ is not the prime character degree graph for any solvable group $G$.
\end{proof}

\end{document}